\documentclass[a4paper, 10pt]{extarticle}
\usepackage{amsmath}
\usepackage{amssymb}
\usepackage{amsfonts}
\usepackage{graphicx}
\usepackage{subfigure}
\usepackage{authblk}
\usepackage{url}
\usepackage[textwidth=16cm,textheight=20cm]{geometry}
\usepackage{hyperref}
\hypersetup{
    colorlinks=true,
    citecolor=blue,
    linkcolor=blue,
    filecolor=blue,
    urlcolor=blue,
}
\newtheorem{theorem}{Theorem}

\newtheorem{remark}[theorem]{Remark}

\title{Singularities in one-dimensional Euler flows}
\author[1\thanks{\textit{E-mail: }\texttt{mihail\underline{ }roop@mail.ru}}]{Mikhail Roop}
\affil[1]{V.A. Trapeznikov Institute of Control Sciences, Russian Academy of Sciences, 65 Profsoyuznaya Str., 117997 Moscow, Russia}

\begin{document}
\maketitle

\abstract{
In this paper, a system of one-dimensional gas dynamics equations is considered. This system is a particular case of Jacobi type systems and has a natural representation in terms of 2-forms on 0-jet space. We use this observation to find a new class of multivalued solutions for an arbitrary thermodynamic state model and discuss singularities of their projections to the space of independent variables for the case of an ideal gas. Caustics and discontinuity lines are found.
}

\section{Introduction}
In this paper, we continue studies of critical phenomena, namely shock waves and phase transitions, appearing in solutions to Euler equations describing flows of gases \cite{LR-ljm-shock,LR-sym}. Such effects have always been paid a great deal of attention both because of their mathematical beauty \cite{Arn1,Arn2,Arn3} and practical applications \cite{Zeld}. The tendency of studying such phenomena is remaining nowadays as well, see, for example \cite{Huang}, where the case of Chaplygin gases is considered, \cite{Ros-Tab,Chat}, where the authors discuss weak shock waves, which is actually the case considered in the present paper, and also it is worth mentioning \cite{Polud1,Polud2}, where the influence of turbulence on detonations is investigated. The properties of global solvability for Euler equations and singularities of their solutions were also studied in a series of works \cite{Tun1,Tun2,Tun3,Tun-Bog}.

Our approach to studying and finding singular properties of solutions to nonlinear PDEs is essentially based on a geometrical theory of PDEs \cite{KrVin,KVL,Ovs,Olver}. Namely, it is known that Euler equations considered here are a particular case of Jacobi type systems (see, for example,\cite{KLR}), which have a natural representation in terms of differential 2-forms on 0-jet space. This observation goes back to a seminal paper \cite{Lych}. One of advantages of this approach is that we need to deal with geometrical structures on 0-jet space instead of 1-jet space, where the equations in question have a natural representation. This idea has also found applications in incompressible hydrodynamics \cite{RoulRub1,RoulRub2}.

Also this approach allows to extend the notion of a smooth solution to a generalized solution understood as an integral manifold of the mentioned forms, which makes it possible to find solutions that are not globally given by functions. Then, singularities of projections of multivalued solutions to the space of independent variables is exactly what corresponds to a formation of shock waves \cite{LychSing}. This concept has been used to describe shock waves in non-stationary filtration problems \cite{AKL-dan,AKL-ifac,AKL-gsa}.

 We combine this observation with another approach to finding multivalued solutions to nonlinear PDEs, which is adding a differential constraint compatible with the original system \cite{KrugLych} that appeared to be fruitful, in particular, in applications to 2-dimensional gas dynamics \cite{Lych-Yum-1}, the Khokhlov-Zabolotskaya equation \cite{Lych-Yum-2}, and also the Hunter-Saxton equation \cite{Eiv}. In the present paper, compared with \cite{LR-ljm-shock,LR-sym}, we give a more accurate description of finding such constraints and discuss two methods of finding them. One of them is based on the mentioned specific geometry that Euler equations have, and another one is based on the theory of differential invariants and quotient PDEs \cite{KLych,Eiv,DLT-sym}. It is worth saying that the last one is more general, however, we decided to discuss both of them to give a more complete picture of the geometry underlying non-stationary Euler equations in one spatial dimension.

The paper is organized as follows. First, we briefly discuss thermodynamics in terms of contact and symplectic geometries \cite{Gibbs,Mrug,Rup,Lych-therm,DLT1,LR-ljm-1}, since it is significant in description of flows of gases. Then, we turn to Euler equations and describe two methods of finding multivalued solutions, using their specific geometric description and using the concept of a quotient PDE. Finally, we get a new class of exact multivalued solutions for an arbitrary thermodynamic model and illustrate them on flows of ideal gases. We find caustics and shock wave fronts. It is worth saying that one can elaborate these solutions for various thermodynamic models, for instance, the van der Waals model taking into account phase transitions, which will be addressed in future papers.

The system of PDEs describing the motion of gases consists of three conservation laws \cite{Batch,Lan-Lif}.

\begin{itemize}
\item Conservation of momentum

\begin{equation}
\label{Euler}
\rho(u_{t}+uu_{x})=-p_{x},
\end{equation}

\item Conservation of mass

\begin{equation}
\label{mass}
\rho_{t}+(\rho u)_{x}=0,
\end{equation}

\item Conservation of entropy along the flow

\begin{equation}
\label{energy}
s_{t}+us_{x}=0,
\end{equation}

\end{itemize}
where $u(t,x)$ is the velocity, $\rho(t,x)$ is the density, $p(t,x)$ is the pressure, $s(t,x)$ is the specific entropy.

One can see that system (\ref{Euler})-(\ref{energy}) is incomplete, which is natural, since the very medium we describe has not yet been specified. This is usually done either by incompressibility condition, when the density is put constant, or by \textit{state equations}, representing the relations between various thermodynamic quantities. Here, we concentrate on the last approach.

\section{Thermodynamics}
We start with the discussion of equilibrium thermodynamics of gases. Any thermodynamical system in equilibrium is described by the following thermodynamical quantities: $e$ --- specific inner energy, $v=\rho^{-1}$ --- specific volume (or inverse of density), $T$ --- temperature, $p$ --- pressure, and $s$ --- specific entropy. The main law of thermodynamics, which is the energy conservation law, states that the differential 1-form
\begin{equation}
\label{theta}
\theta=-ds+T^{-1}de-pT^{-1}\rho^{-2}d\rho
\end{equation}
must vanish. This drives us to the notion of \textit{thermodynamic states} understood as maximal integral manifolds of \eqref{theta}.

More precisely, consider the contact manifold $\widehat{\Phi}=W\times W^{\ast}\times\mathbb{R}\simeq\mathbb{R}^{5}$, where $W\simeq\mathbb{R}^{2}(e,\rho)$, $W^{\ast}\simeq\mathbb{R}^{2}(p,T)$, and $s\in\mathbb{R}$. Then, a thermodynamic state is a Legendrian manifold $\widehat{L}\subset(\widehat{\Phi},\theta)$, such that
\begin{equation*}
\theta|_{\widehat{L}}=0.
\end{equation*}
This exactly means that the energy conservation law holds on $\widehat{L}$.

Let us choose $(e,\rho)$ as local coordinates on $\widehat{L}$. Then, the two-dimensional manifold $\widehat{L}\subset(\widehat{\Phi},\theta)$ is given by three relations
\begin{equation}
\label{legmani}
\widehat{L}=\left\{s=S(e,\rho),\, T=\frac{1}{S_{e}},\, p=-\rho^{2}\frac{S_{\rho}}{S_{e}}\right\}
\end{equation}
for some function $S(e,\rho)$.

However, in practice there are no ways of determining $S(e,\rho)$, which motivates us to switch to the Lagrangian viewpoint. Namely, consider projection
\begin{equation*}
\pi\colon\widehat{\Phi}\to\Phi,\quad\pi(p,T,e,\rho,s)=(p,T,e,\rho).
\end{equation*}
A pair $(\Phi,\Omega)$ is a symplectic manifold with the structure form
\begin{equation*}
\Omega=d\theta=d(T^{-1})\wedge de-d(pT^{-1}\rho^{-2})\wedge d\rho.
\end{equation*}

Then, a thermodynamic state is a Lagrangian manifold $\pi(\widehat{L})=L\subset(\Phi,\Omega)$, such that
\begin{equation*}
\Omega|_{L}=0.
\end{equation*}
In the symplectic space $(\Phi,\Omega)$ the Lagrangian manifold $L$ is given by \textit{state equations}:
\begin{equation*}L=\left\{f(p,T,e,\rho)=0,\,g(p,T,e,\rho)=0\right\},\end{equation*}
such that the Poisson bracket $[f,g]$ with respect to the structure form $\Omega$
\begin{equation*}[f,g]\,\Omega\wedge\Omega=df\wedge dg\wedge\Omega\end{equation*}
vanishes on $L$:
\begin{equation}
\label{compat}
[f,g]=0\text{ on }L.
\end{equation}
Condition $[f,g]=0$ on $L$ is called the \textit{compatibility condition} for state equations.

If one chooses $(T,\rho)$ as local coordinates on the Lagrangian manifold $L$, that is
\begin{equation*}
L=\left\{p=P(T,\rho),\,e=E(T,\rho)\right\}.
\end{equation*}
The condition $[f,g]=0$ on $L$ leads to the equation
\begin{equation*}(-\rho^{-2}T^{-1}P)_{T}=(T^{-2}E)_{\rho},
\end{equation*}
and therefore the following theorem is valid:
\begin{theorem}
The Lagrangian manifold $L$ is given by means of \textit{the Massieu-Planck potential} $\phi(\rho,T)$
\begin{equation}
\label{lag-mani}
p=-\rho^{2}T\phi_{\rho},\quad e=T^{2}\phi_{T}.
\end{equation}
\end{theorem}
\begin{remark}
One can build up a Legendrian manifold \eqref{legmani} from a given Lagrangian one by resolving the overdetermined system of equations on function $S(e,\rho)$:
\begin{equation}
\label{sys-ent}
T=\frac{1}{S_{e}},\quad p=-\rho^{2}\frac{S_{\rho}}{S_{e}},
\end{equation}
where $p(e,\rho)$ and $T(e,\rho)$ are specified once \eqref{lag-mani} is given, and the compatibility condition for \eqref{sys-ent} is \eqref{compat}. Indeed, the specific entropy is expressed in terms of the Massieu-Planck potential as follows \cite{LR-ljm-1}:
\begin{equation*}
s=\phi+T\phi_{T}.
\end{equation*}
\end{remark}

The manifold $\Phi=W\times W^{\ast}$ is equipped with the differential quadratic form \cite{Lych-therm}
\begin{equation*}
\kappa=d(T^{-1})\cdot de-\rho^{-2}d(pT^{-1})\cdot d\rho.
\end{equation*}
Not all the points on the Lagrangian manifold correspond to real physical states, but only those where (applicability condition)
\begin{equation*}
\kappa|_{L}<0.
\end{equation*}
Writing down $\kappa|_{L}$ in terms of the Massieu-Planck potential, we get:
\begin{equation*}
\kappa|_{L}=-(2T^{-1}\phi_{T}+\phi_{TT})dT\cdot dT+(2\rho^{-1}\phi_{\rho}+\phi_{\rho\rho})d\rho\cdot d\rho.
\end{equation*}
Taking into account \eqref{lag-mani} we observe that the applicability condition $\kappa|_{L}<0$ is
\begin{equation*}
e_{T}>0,\quad p_{\rho}>0,
\end{equation*}
which is known to be the conditions of thermodynamic stability.

In the context of the geometrical approach to thermodynamics, \textit{thermodynamic processes} are understood as contact transformations of $\widehat{\Phi}$, preserving the Legendrian manifold $\widehat{L}$. They are generated by contact vector fields, and their integral curves $l\subset\widehat{L}$ we will also call \textit{thermodynamic processes}.
\section{Euler equations}
Let us note that system (\ref{Euler})-(\ref{energy}) becomes complete once we fix a thermodynamic state \eqref{legmani}.

We will apply a \textit{homentropicity ansatz}, which means that we put $s(t,x)=s_{0}=\mathrm{const}$. On the one hand, we get rid of equation \eqref{energy}, on the other hand, we are able to express all the thermodynamic quantities in terms of the density $\rho$. Indeed, considering $s_{0}=\phi+T\phi_{T}$ as an equation for $T(\rho)$, we observe that the derivative of its right-hand side with respect to $T$ is positive in an applicable domain, and therefore due to the implicit function theorem, this equation determines $T(\rho)$ uniquely. By means of state equations we also get $p=p(\rho)$. Summarizing above discussion, we end up with the following two-component system of PDEs:

\begin{equation}
\label{EU}
\mathcal{E}=\begin{cases}
 \rho_{t}+(\rho u)_{x}=0,
   \\
\displaystyle u_{t}+uu_{x}+\frac{p^{\prime}(\rho)}{\rho}\rho_{x}=0.
 \end{cases}
\end{equation}
\subsection{Geometrical structures associated to Euler equations}
Let $E=J^{0}(t,x,u,\rho)$ be the space of 0-jets, and $M=\mathbb{R}^{2}(t,x)$. Then, following \cite{KLR,Lych}, one can associate two differential 2-forms
\begin{equation*}
\begin{cases}
 \omega_{1}=\rho dt\wedge du+udt\wedge d\rho-dx\wedge d\rho,
   \\
\displaystyle \omega_{2}=udt\wedge du+\frac{p^{\prime}(\rho)}{\rho}dt\wedge d\rho-dx\wedge du
 \end{cases}
\end{equation*}
with system \eqref{EU}. Indeed, let $\Omega^{2}(E)$ be a module of differential 2-forms on $E$, then any differential 2-form $\omega\in\Omega^{2}(E)$ generates an operator
\begin{equation*}
\Delta_{\omega}\colon C^{\infty}(M)\to\Omega^{2}(M),\quad \Delta_{\omega}(f)=\omega|_{\Gamma^{0}(f)},
\end{equation*}
where $\Gamma^{0}(f)\subset E$ is a graph of the vector-function $f$. The system $\mathcal{E}$ can now be written as
\begin{equation*}
\Delta_{\omega_{1}}(f)=0,\quad\Delta_{\omega_{2}}(f)=0,
\end{equation*}
where $f=(u(t,x),\rho(t,x))$.

A 2-dimensional manifold $N\subset E$ is said to be a \textit{(multivalued) solution} of $\mathcal{E}$, if $\omega_{1}|_{N}=\omega_{2}|_{N}=0$.

Note that the form $\omega_{2}$ is closed and non-degenerate, and therefore may serve as a symplectic structure on $E$, which means that any multivalued solution to $\mathcal{E}$ is a Lagrangian manifold.

Let us fix a volume form on $E$ as $q=dt\wedge dx\wedge du\wedge d\rho$ and introduce a bilinear operator
\begin{equation*}
P\colon\Omega^{2}(E)\times\Omega^{2}(E)\to C^{\infty}(E),\quad \alpha_{1}\wedge\alpha_{2}=P(\alpha_{1},\alpha_{2})q,\quad\alpha_{i}\in\Omega^{2}(E).
\end{equation*}
Let $P_{\omega}=\|P(\omega_{i},\omega_{j})\|$ be its matrix in the basis $\langle\omega_{1},\omega_{2}\rangle$. Then, system $\mathcal{E}$ is said to be \textit{hyperbolic}, if $\det(P_{\omega})<0$, \textit{elliptic}, if $\det(P_{\omega})>0$ and \textit{parabolic}, if $\det(P_{\omega})=0$. In our case the matrix $P_{\omega}$ has the following form:
\begin{equation*}
P_{\omega}=
\begin{pmatrix}
2\rho & 0\\
0 & -2\rho^{-1}p^{\prime}(\rho)
\end{pmatrix},
\end{equation*}
and therefore the condition for the system $\mathcal{E}$ to be of hyperbolic type is $p^{\prime}(\rho)>0$.

It is worth saying that the applicability condition $e_{T}>0$ is satisfied for a considerable number of real gas models on the entire Lagrangian manifold, and we will consider only such models, while another applicability condition $p_{\rho}>0$ violates in some domains of Lagrangian manifolds of real gases, for instance, van der Waals model \cite{LR-sym}. Thus the applicability condition $p_{\rho}>0$ is essential for us both in the context of thermodynamics and in the context of type of the Euler system.
\begin{theorem}
If a thermodynamic process curve $l\subset\widehat{L}$ lies in the domain of negativity of the form $\kappa|_{L}$, then the system $\mathcal{E}$ is hyperbolic.
\end{theorem}

Obviously, the forms $\omega_{1}$ and $\omega_{2}$ are defined up to a non-degenerate linear transformation. Indeed, the forms
\begin{eqnarray*}
\label{forms}
\widehat{\omega_{1}}&=&a_{11}\omega_{1}+a_{12}\omega_{2},\\
\widehat{\omega_{2}}&=&a_{21}\omega_{1}+a_{22}\omega_{2},
\end{eqnarray*}
where $a_{11}a_{22}-a_{12}a_{21}\ne 0$, define the Euler system as well. In a hyperbolic case one can choose $a_{ij}$ in such a way that
\begin{equation}
\label{eff}
\omega_{1}\wedge\omega_{2}=0,\quad \omega_{1}\wedge\omega_{1}=-\omega_{2}\wedge\omega_{2}.
\end{equation}
Differential 2-forms satisfying relations \eqref{eff} are also called \textit{effective} \cite{KLR}. For the case of Euler system, we get the following theorem \cite{LR-ljm-shock}:
\begin{theorem}
Let $\mathcal{E}$ be a system of hyperbolic type. Then, it can be given by 2-forms
\begin{eqnarray*}
\omega_{1}&=&A(\rho)(\rho dt\wedge du+udt\wedge d\rho-dx\wedge d\rho),\\
\omega_{2}&=&udt\wedge du+\rho A^{2}(\rho)dt\wedge d\rho-dx\wedge du,
\end{eqnarray*}
where $A(\rho)=\rho^{-1}\sqrt{p^{\prime}(\rho)}$, and 2-forms $\omega_{1}$, $\omega_{2}$ satisfy relations (\ref{eff}).
\end{theorem}

Let us introduce another operator
\begin{equation*}
A_{\omega}\colon D(E)\to D(E),\quad X\rfloor\omega_{2}=A_{\omega}(X)\rfloor\omega_{1},
\end{equation*}
where $D(E)$ is the module of vector fields on $E$, and $X\in D(E)$.

If one chooses $\langle\partial_{t},\partial_{x},\partial_{\rho},\partial_{u}\rangle$ as a basis in $D(E)$, one gets the matrix of this operator
\begin{equation*}
W=\frac{1}{\rho A(\rho)}\begin{pmatrix}
u & -1 & 0 & 0\\
u^{2}-\rho^{2}A^{2}(\rho) & -u & 0 & 0\\
0 &0 & 0 & \rho A^{2}(\rho)\\
0 &0 & \rho & 0
\end{pmatrix},
\end{equation*}
and moreover, in a hyperbolic case $A_{\omega}^{2}=\mathrm{id}$.

Let $a\in E$, $\mathcal{C}_{\pm}(a)$ be eigenspaces of the operator $A_{\omega}(a)$, then
\begin{equation*}
T_{a}E=\mathcal{C}_{-}(a)\oplus\mathcal{C}_{+}(a),
\end{equation*}
and \textit{characteristic distributions} $\mathcal{C}_{+}=\langle X_{+},Y_{+}\rangle$ and $\mathcal{C}_{-}=\langle X_{-},Y_{-}\rangle$ are generated by vector fields \cite{LR-ljm-shock}
\begin{eqnarray*}
X_{\pm}&=&\pm A(\rho)\partial_{u}+\partial_{\rho},\\
Y_{\pm}&=&(\mp \rho A(\rho)+u)^{-1}\partial_{t}+\partial_{x}.
\end{eqnarray*}
In the case when both distributions $\mathcal{C}_{+}$ and $\mathcal{C}_{-}$ are integrable, one can explicitly solve the Cauchy problem (see, for example, \cite{KLR}). The integrability conditions for $\mathcal{C}_{+}$ and $\mathcal{C}_{-}$ are given by the following theorem \cite{LR-ljm-shock}:
\begin{theorem}
Distributions $\mathcal{C}_{+}$ and $\mathcal{C}_{-}$ are integrable if
\begin{equation}
\label{intr}
p(\rho)=c_{0}\rho^{3}+c_{1},
\end{equation}
where $c_{0}$ and $c_{1}$ are constants.
\end{theorem}

\begin{figure}[h!]
\centering
\includegraphics[scale=.5]{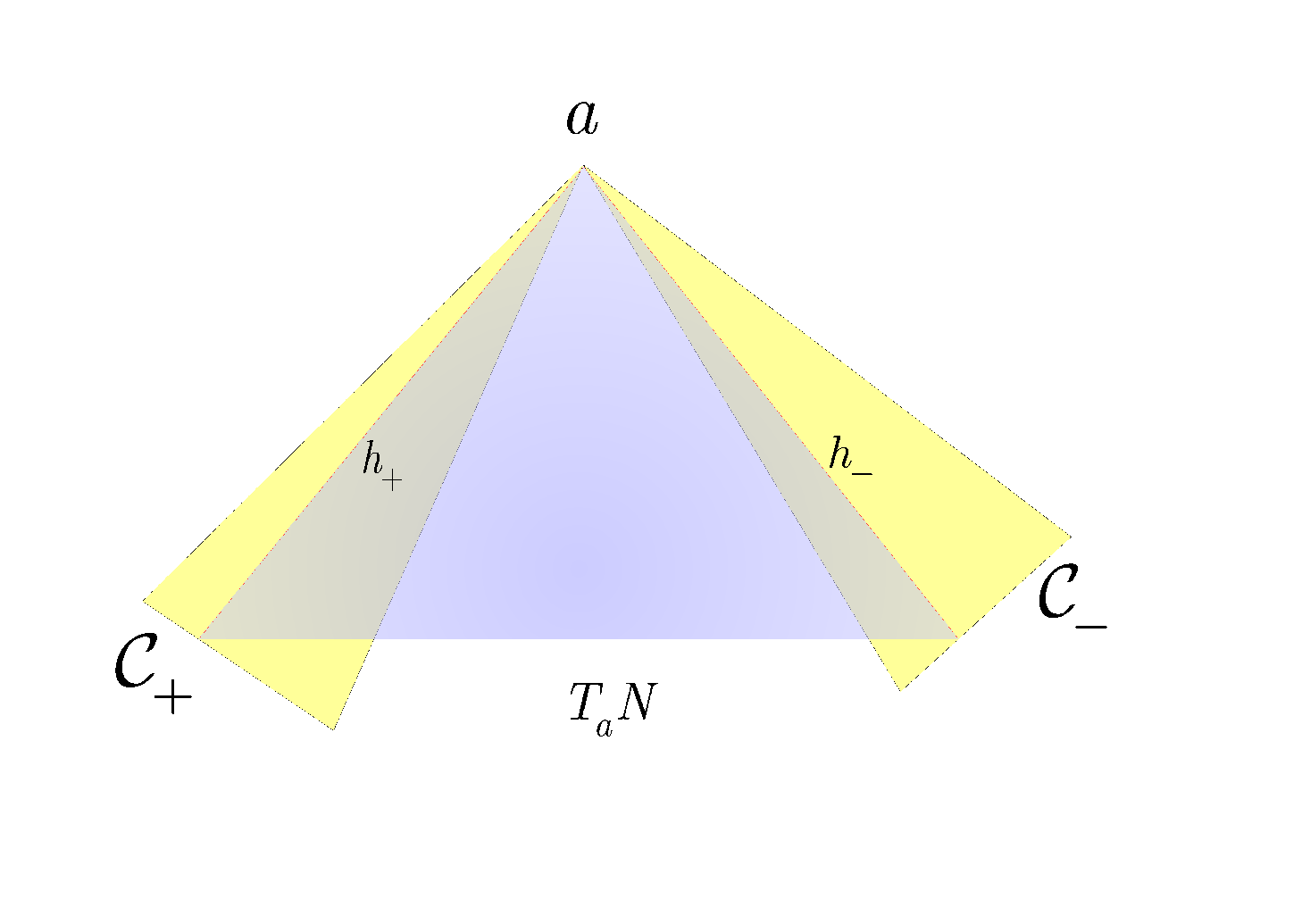}
\caption{Characteristic directions.}
\label{char-dir}
\end{figure}

Nevertheless, we will not specify the dependence $p(\rho)$. To construct solutions in general case, we will need the following theorem \cite{KLR}:
\begin{theorem}
\label{thm-1}
A two-dimensional manifold $N\subset E$ is a multivalued solution to $\mathcal{E}$, if and only if the tangent spaces $T_{a}N$ for all $a\in N$ have one-dimensional intersections $h_{+}(a)$ and $h_{-}(a)$ with planes $\mathcal{C}_{+}(a)$ and $\mathcal{C}_{-}(a)$:
\begin{equation*}
T_{a}N=h_{+}(a)\oplus h_{-}(a).
\end{equation*}
\end{theorem}
The lines $h_{\pm}(a)=\mathcal{C}_{\pm}(a)\cap T_{a}N$ are called \textit{characteristic directions}. The statement of the above theorem is graphically shown in Fig.~\ref{char-dir}.

We will look for a one-parametric family of solutions defined by some 3-dimensional manifold
\begin{equation*}
N_{0}=\left\{F(t,x,u,\rho)=0\right\}
\end{equation*}
with integrable distribution $h=\langle h_{+},h_{-}\rangle$.

Let us take two vector fields $V_{\pm}\in\mathcal{C}_{\pm}$, tangent to $N_{0}$, i.e. $V_{\pm}\rfloor dF=0$ and choose $F(t,x,u,\rho)$ in such a way that $V_{+}\wedge V_{+}\wedge[V_{+},V_{-}]=0$. Restricting $V_{\pm}$ the $N_{0}$, we get an integrable distribution $h=\langle h_{+},h_{-}\rangle$, where $h_{+}=V_{+}|_{N_{0}}$ and $h_{-}=V_{-}|_{N_{0}}$.

Consider a particular case $F(t,x,\rho,u)=t-f(\rho,u)/\rho$. Then the integrability condition $V_{-}\wedge V_{+}\wedge[V_{-},V_{+}]=0$ leads to the following equation for $f(\rho,u)$:
\begin{equation}
\label{wave}
f_{uu}-A^{-2}(\rho)f_{\rho\rho}=0,
\end{equation}
which in some cases can be reduced to the wave equation with constant coefficients.
\begin{theorem}
Equation (\ref{wave}) is equivalent to the wave equation with constant coefficients if
\begin{equation*}
A(\rho)=(\beta_{1}\rho+\beta_{2})^{-2},
\end{equation*}
where $\beta_{1}$, $\beta_{2}$ are constants.
\end{theorem}
In this case we are able to find all 3-dimensional manifolds $N_{0}$ of the form $N_{0}=\left\{t-f(\rho,u)/\rho=0\right\}$.

Taking any of solutions to \eqref{wave} we get the corresponding integrable distribution $h=\langle h_{+},h_{-}\rangle$, and according to Theorem~\ref{thm-1}, its integrals give us multivalued solutions.
\subsection{Integrability via quotients}
There is another interpretation of \eqref{wave} as a quotient PDE for the system $\mathcal{E}$. To show this, let us first recall some concepts from the geometrical theory of PDEs following \cite{KLych,KrVin,KVL}.

Let $(x_{1},\ldots,x_{n})$ be a collection of independent variables and let $(u^{1},\ldots,u^{m})$ be a collection of dependent variables. Let $\mathcal{E}\subset J^{k}(\mathbb{R}^{n})$ be a system of PDEs with the symmetry Lie algebra $\mathfrak{g}$. By $\mathcal{E}_{l}\subset J^{l}(\mathbb{R}^{n})$ for some $l\ge k$ we will denote the prolongation of $\mathcal{E}$. Then, under some conditions, the field of rational differential $\mathfrak{g}$-invariants is finitely generated. More precisely, the global Lie-Tresse theorem \cite{KLych} is valid.
\begin{theorem}[Kruglikov, Lychagin]
Let $\mathcal{E}_{k}\subset J^{k}(\mathbb{R}^{n})$ be an algebraic formally integrable differential equation and let $\mathfrak{g}$ be its algebraic symmetry Lie algebra. Then, there exist rational differential $\mathfrak{g}$-invariants $a_{1},\ldots,a_{n},b^{1},\ldots,b^{N}$ of order $\le l$, such that the field of rational $\mathfrak{g}$-invariants is generated by rational functions of these functions and Tresse derivatives $\frac{d^{|\alpha|}b^{j}}{da^{\alpha}}$.
\end{theorem}
Remark that conditions of the global Lie-Tresse theorem are not very restrictive, and its statement is true for a considerable number of PDE systems and their symmetry Lie algebras, in particular, for Euler equations.

The Tresse derivatives mentioned in the Lie-Tresse theorem are constructed as follows (we omit here some technical details while emphasizing on general concepts, and we refer to \cite{KLych} for details). Let us take $n$ horizontally independent differential invariants $a_{1},\ldots,a_{n}$, which means that
\begin{equation}
\label{indep-cond}
\widehat{d}a_{1}\wedge\ldots\wedge\widehat{d}a_{n}\ne0
\end{equation}
in some Zariski-open set in $J^{l+1}(\mathbb{R}^{n})$. In \eqref{indep-cond}, $\widehat{d}$ is the total differential:
\begin{equation*}
\widehat{d}f=\sum\limits_{i=1}^{n}\frac{df}{dx_{i}}dx_{i},
\end{equation*}
where
\begin{equation*}
\frac{d}{dx_{i}}=\frac{\partial}{\partial x_{i}}+\sum\limits_{j,\sigma}u^{j}_{\sigma i}\frac{\partial}{\partial u_{\sigma}^{j}}.
\end{equation*}
Tresse derivatives are constructed as partial derivatives with respect to invariants $a_{1},\ldots,a_{n}$:
\begin{equation*}
\frac{d}{da_{i}}=\sum\limits_{j}A_{ij}\frac{d}{dx_{j}},\quad \frac{da_{j}}{da_{i}}=\delta_{ij},
\end{equation*}
where $\delta_{ij}$ is the Kronecker delta. Condition \eqref{indep-cond} guarantees the existence of solution to the system of linear equations on $A_{ij}$. By applying the Tresse derivatives to differential invariants we get new differential invariants.

In general, the algebra of differential invariants is not freely generated, there are relations between invariants, called \textit{syzygies}.

Let $s$ be a solution to $\mathcal{E}$, and consider restrictions:
\begin{equation*}
a_{1}|_{s},\ldots,a_{n}|_{s},b^{1}|_{s},\ldots,b^{N}|_{s},
\end{equation*}
which locally can be viewed as $(n+N)$ functions on an $n$-dimensional manifold, therefore (locally)
\begin{equation*}
b^{j}|_{s}=B^{j}(a_{1}|_{s},\ldots,a_{n}|_{s}),\quad j=1,\ldots,N.
\end{equation*}
In fact, $B^{j}$ depend on the equivalence class of $s$ (where the equivalence relation is defined by the symmetry Lie algebra), rather than on $s$ itself.

Removing restrictions, we get:
\begin{equation}
\label{rel}
b^{j}=B^{j}(a_{1},\ldots,a_{n}).
\end{equation}
Functions $B^{j}$ can be found from \textit{quotient PDEs}. Let us apply the Tresse derivatives to \eqref{rel}. We get
\begin{equation*}
J_{ji}=\frac{db^{j}}{da_{i}}=B^{j}_{a_{i}}(a_{1},\ldots,a_{n}),\ldots,J_{j\alpha\beta}=\frac{d^{2}b^{j}}{da_{\alpha}da_{\beta}}=B^{j}_{a_{\alpha}a_{\beta}}(a_{1},\ldots,a_{n}),\ldots
\end{equation*}
Since $b^{j}$ are invariants, $\frac{d}{da_{i}}$ are invariant derivations, functions $J_{ji},J_{j\alpha\beta},\ldots$ are invariants too. Finding syzygies between invariants $a_{i}$, $b^{j}$, $J_{ji}$, $J_{j\alpha\beta},\ldots$ we get a relation (perhaps, a number of them)
\begin{equation*}
F(a_{1},\ldots,a_{n},B^{1},\ldots,B^{N},B^{j}_{a_{i}},B^{j}_{a_{\alpha}a_{\beta}},\ldots)=0,
\end{equation*}
called a \textit{quotient PDE}, which is a PDE on functions $B^{j}(a_{1},\ldots,a_{n})$.

Let us collect some of the most important properties of quotient PDEs:
\begin{enumerate}
\item Solution to a quotient PDE is a $\mathfrak{g}$-orbit of $s$;
\item Solutions to quotient PDEs provide us with differential constraints $b^{j}-B^{j}(a_{1},\ldots,a_{n})=0$, compatible with the original system $\mathcal{E}$.
\end{enumerate}
Finding such constraints, we reduce the integration of $\mathcal{E}$ to the integration of a completely integrable Cartan distribution $\mathcal{C}_{\mathcal{E}}$ with the same symmetry algebra.

Let us now apply these ideas to the Euler system $\mathcal{E}$.
\begin{theorem}
The symmetry Lie algebra of the system of Euler equations \eqref{EU} is generated by vector fields
\begin{equation*}
\begin{split}
X_{1}=\frac{\partial}{\partial x},\quad &X_{2}=\frac{\partial}{\partial t},\quad X_{3}=t\frac{\partial}{\partial x}+\frac{\partial}{\partial u},\quad X_{4}=t\frac{\partial}{\partial t}+x\frac{\partial}{\partial x},\quad {}\\&
X_{5}=C_{1}e^{-C_{1}u}Y\frac{\partial}{\partial t}+(C_{1}Yu+Y^{\prime}\rho+Y)e^{-C_{1}u}\frac{\partial}{\partial x},
\end{split}
\end{equation*}
where $Y=Y(\rho)$ is a solution to the ODE
\begin{equation*}
-\rho A^{2}(\rho)C_{1}^{2}Y+\rho Y^{\prime\prime}+2Y^{\prime}=0,
\end{equation*}
and $C_{1}$ is a constant.
\end{theorem}
Let us choose $\mathfrak{g}=\langle\partial_{x}\rangle$. Then, we have the following 0-order invariants:
\begin{equation*}
a_{1}=u,\quad a_{2}=\rho,\quad b^{1}=t.
\end{equation*}
Put $b^{1}=t=B(u,\rho)$. Then, Tresse derivatives will be:
\begin{eqnarray*}
\frac{d}{du}&=&\frac{\rho_{x}}{\rho u_{x}^{2}-\rho A^{2}\rho_{x}^{2}}\frac{d}{dt}+\frac{\rho u_{x}+u\rho_{x}}{\rho u_{x}^{2}-\rho A^{2}\rho_{x}^{2}}\frac{d}{dx},\\
\frac{d}{d\rho}&=&\frac{u_{x}}{\rho A^{2}\rho_{x}^{2}-\rho u_{x}^{2}}\frac{d}{dt}+\frac{A^{2}\rho\rho_{x}+uu_{x}}{\rho A^{2}\rho_{x}^{2}-\rho u_{x}^{2}}\frac{d}{dx}.
\end{eqnarray*}
Finding syzygies between $u$, $\rho$, $\frac{dB}{du}$, $\frac{dB}{d\rho}$, $\frac{d^{2}B}{du^{2}}$, $\frac{d^{2}B}{d\rho^{2}}$, $\frac{d^{2}B}{du d\rho}$, we get a quotient PDE:
\begin{equation*}
A^{2}\rho\frac{d^{2}B}{du^{2}}-\rho\frac{d^{2}B}{d\rho^{2}}-2\frac{dB}{d\rho}=0.
\end{equation*}
Putting $B(u,\rho)=f(u,\rho)/\rho$, we get
\begin{equation}
\label{QUO-PDE}
f_{uu}-A^{-2}f_{\rho\rho}=0,
\end{equation}
which coincides with \eqref{wave}.

Let us take any solution $f(u,\rho)$ to \eqref{QUO-PDE}. Then, differentiating $t-f(u,\rho)/\rho=0$ with respect to $t$ and $x$, we get two more PDEs, and together with Euler equations we get 5 relations on 1-jet space, determining the submanifold $\widehat{\mathcal{E}}\subset J^{1}(\mathbb{R}^{2})$.

So, $\dim J^{1}(\mathbb{R}^{2})=8$, $\dim\widehat{\mathcal{E}}=3$, $\dim\mathcal{C}_{\widehat{\mathcal{E}}}=2\Longrightarrow\mathrm{codim}\,\mathcal{C}_{\widehat{\mathcal{E}}}=1$, and the Cartan distribution $\mathcal{C}_{\widehat{\mathcal{E}}}$ is given by a differential 1-form
\begin{equation}
\label{Cart-form}
\begin{split}
\omega=&-\frac{\rho  (f_{\rho} \rho -f)}{\rho^2 A^2 f_{u}^2-f_{\rho}^2 \rho^2+2 f_{\rho} f\rho-f^2}dx+\frac{(f_{\rho} \rho -f) (f_{u} u -f_{\rho} \rho +f)}{\rho^2 A^2 f_{u}^2-f_{\rho}^2 \rho ^2+2 f_{\rho} f \rho -f^2}du+{}\\&+\frac{(-\rho ^2 A^2 f_{u}+f_{\rho} \rho  u -f u ) (f_{\rho} \rho -f)}{\rho^3 A^2 f_{u}^2-\rho^3 f_{\rho}^2+2 \rho^2 f f_{\rho}-f^2 \rho}d\rho,
\end{split}
\end{equation}
such that $\omega\wedge d\omega=0$ for any solution $f(u,\rho)$ to \eqref{wave}. Integrals of the form $\omega$ give us solutions to $\widehat{\mathcal{E}}$.
\begin{remark}
The distribution generated by $\omega$ coincides with that generated by $\langle h_{-},h_{+}\rangle$.
\end{remark}
\subsection{Solutions}
We can see that finding differential constraint to integrate the Euler system can be performed in one of two equivalent ways described in this paper. We will proceed with the integration of $\mathcal{E}$ by finding integrals of the form \eqref{Cart-form}.

Note that the vector field $\partial_{x}$ is an infinitesimal transversal symmetry of \eqref{Cart-form}, and therefore the differential form $\varkappa=\omega/\omega(\partial_{x})$ is closed and therefore locally exact, i.e. $\varkappa=dJ$ for some function $J\in C^{\infty}(E)$, and $J=\mathrm{const}$ together with $t-f(u,\rho)/\rho=0$ give us solutions to $\mathcal{E}$.

Let us apply the following separation ansatz to \eqref{wave}: $f(u,\rho)=\rho\mu(u)+\nu(\rho)$. Then, we get two ODEs on $\mu(u)$ and $\nu(\rho)$:
\begin{equation*}
\mu^{\prime\prime}=\lambda,\quad \nu^{\prime\prime}=\lambda\rho A^{2},
\end{equation*}
for some constant $\lambda$, and we get the first quadrature for $\mathcal{E}$:
\begin{equation}
\label{sol-t}
t=\frac{\lambda u^{2}}{2}+\alpha_{0}u+\frac{\alpha_{2}}{\rho}+\frac{\lambda}{\rho}\int Q(\rho)d\rho+t_{0},
\end{equation}
where $\alpha_{j}, t_{0}$ are constants,  $Q(\rho)=\int\rho A^{2}d\rho$.

The differential 1-form $\varkappa$ equals
\begin{equation}
\label{varkappa}
\begin{split}
\varkappa=&dx-\left(\lambda u^{2}+\alpha_{0}u-\lambda Q(\rho)+\frac{\lambda}{\rho}\int Q(\rho)d\rho+\frac{\alpha_{2}}{\rho}\right)du+{}\\&+\left(\frac{\lambda u}{\rho^{2}}\int Q(\rho)d\rho-\frac{\lambda u}{\rho}Q(\rho)+\rho A^{2}(\rho)(\lambda u+\alpha_{0})+\frac{\alpha_{2}u}{\rho^{2}}\right)d\rho.
\end{split}
\end{equation}
Integrating \eqref{varkappa}, we get the second quadrature for $\mathcal{E}$:
\begin{equation}
\label{sol-x}
x=\frac{\lambda u^{3}}{3}-\lambda uQ(\rho)+\frac{\alpha_{0}u^{2}}{2}+\frac{\lambda u}{\rho}\int Q(\rho)d\rho+\frac{\alpha_{2}u}{\rho}-\alpha_{0}Q(\rho)+x_{0},
\end{equation}
where $x_{0}$ is a constant.

Formulae \eqref{sol-t},\eqref{sol-x} give us a multivalued solution to $\mathcal{E}$. Note that it can be applied to any thermodynamic state model, because $p(\rho)$ is not specified.

Singularities of projections of the multivalued solution $N$ given by \eqref{sol-t},\eqref{sol-x} to the space of independent variables appear where the differential 2-form $dx\wedge dt$ degenerates. Equation $(dx\wedge dt)|_{N}=0$ gives us a curve, called \textit{caustic}. Parametric equations for the caustic are
\begin{eqnarray*}
x(\rho)&=&x_{0}-\frac{Z_{\pm}^{3}}{3\lambda^{2}A^{3}\rho^{6}}+\frac{\alpha_{0}Z_{\pm}^{2}}{2\lambda^{2}A^{2}\rho^{4}}\pm\frac{\alpha_{0}Z_{\pm}}{\lambda\rho}\mp\frac{Z_{\pm}^{2}}{\lambda A\rho^{3}}-\alpha_{0}Q,\\
t(\rho)&=&\frac{Z_{+}^{2}}{2\lambda A^{2}\rho^{4}}+\frac{Z_{+}}{\rho}\left(1-\frac{\alpha_{0}}{\lambda\rho A}\right)-\alpha_{0}\rho A+\lambda Q+t_{0},
\end{eqnarray*}
where
\begin{equation*}
Z_{\pm}=\alpha_{0}\rho^{2}A\mp\lambda\rho Q\pm\lambda\int Qd\rho\pm\alpha_{2}.
\end{equation*}

To cut a multivalued solution given by \eqref{sol-t},\eqref{sol-x} into single valued branches representing the corresponding discontinuous solution, one needs a conservation law. Following \cite{Lan-Lif}, we will use the mass conservation law. Let us choose $(t,\rho)$ as coordinates on the solution $N$. Resolving \eqref{sol-t} with respect to $u$, we get two roots
\begin{equation}
\label{u-rho-t-sol}
u=U_{\pm}(\rho,t)=-\frac{\alpha_{0}}{\lambda}\pm\frac{1}{\lambda\rho}\sqrt{2\rho\left(\rho\lambda(t-t_{0})+\frac{\rho\alpha_{0}^{2}}{2}-\lambda\alpha_{2}-\lambda^{2}\int Qd\rho\right)}.
\end{equation}

Substituting this relation to \eqref{mass}, we get
\begin{equation*}
\rho_{t}+(\rho U_{\pm})_{x}=0,
\end{equation*}
and therefore the conservation law is
\begin{equation}
\label{con-law}
\Theta_{\pm}=\rho dx-\rho U_{\pm}dt.
\end{equation}
The restriction $\Theta_{\pm}|_{N}$ of \eqref{con-law} to solution \eqref{sol-t},\eqref{sol-x} is a closed form, and therefore locally $\Theta_{\pm}|_{N}=dH_{\pm}$, and the potential $H_{\pm}(\rho,t)$ equals
\begin{equation*}
H_{\pm}(\rho,t)=\mp\frac{1}{\lambda\sqrt{\rho}}\sqrt{2\left(\rho\lambda(t-t_{0})+\frac{\rho\alpha_{0}^{2}}{2}-\lambda\alpha_{2}-\lambda^{2}\int Qd\rho\right)}\left(\rho\lambda Q-\lambda\int Q d\rho-\alpha_{2}\right).
\end{equation*}

The points where solution has discontinuity (shock wave front) is found from the system of equations:
\begin{equation*}
H_{\pm}(\rho_{1},t)=H_{\pm}(\rho_{2},t), \quad x(\rho_{1},t)=x(\rho_{2},t),
\end{equation*}
which is to be resolved at each time moment $t$ with respect to $\rho_{1,2}$ --- the values of densities on each branch of a discontinuous solution, and where $x(\rho,t)$ is found as a result of substitution of \eqref{u-rho-t-sol} to \eqref{sol-x}.

\section{Solutions for an ideal gas, caustics, shock waves}
Here, we illustrate solutions obtained in the previous section on an ideal gas model. The Legendrian manifold for an ideal gas is given by
\begin{equation*}
\widehat{L}=\left\{p=R\rho T,\, e=\frac{n}{2}RT,\, s=R\ln\left(\frac{T^{n/2}}{\rho}\right)\right\},
\end{equation*}
where $R$ is the universal gas constant, and $n$ is the degree of freedom.

The differential quadratic form $\kappa|_{L}$ for ideal gases is
\begin{equation*}
\kappa|_{L}=-\frac{Rn}{2}\frac{dT^{2}}{T^{2}}-R\rho^{-2}d\rho^{2},
\end{equation*}
and it is negative on the entire Lagrangian manifold. Therefore, the system $\mathcal{E}$ in case of ideal gases is hyperbolic for any process $l\subset\widehat{L}$.

By fixing the entropy constant $s=s_{0}$, we get the following expressions for $T(\rho)$ and $p(\rho)$:
\begin{equation*}
T(\rho)=\exp\left(\frac{2s_{0}}{Rn}\right)\rho^{2/n},\quad p(\rho)=R\exp\left(\frac{2s_{0}}{Rn}\right)\rho^{2/n+1},
\end{equation*}
and the function $A(\rho)=\rho^{-1}\sqrt{p^{\prime}(\rho)}=A_{0}\rho^{m}$ (see also \cite{LR-ljm-shock}), where
\begin{equation*}
A_{0}=\sqrt{R\left(1+\frac{2}{n}\right)\exp\left(\frac{2s_{0}}{Rn}\right)},\quad m=\frac{1}{n}-1.
\end{equation*}

In the case of an ideal gas, we get the following quadratures:
\begin{equation*}
t=\frac{\lambda u^{2}}{2}+\alpha_{0}u+\frac{\alpha_{2}}{\rho}+\frac{\lambda A_{0}^{2}\rho^{2m+2}}{2(2m^{2}+5m+3)}+t_{0},
\end{equation*}
\begin{equation*}
\begin{split}
x=&\frac{1}{\rho(2m^{2}+5m+3)}\left((m+1)(2m+3)\left(\alpha_{2}u+\rho\left(\frac{\lambda u^{3}}{3}+\frac{\alpha_{0}u^{2}}{2}+x_{0}\right)\right)-\right.{}\\&-\left.A_{0}^{2}\rho^{2m+3}\left(m(\lambda u+\alpha_{0})+\lambda u+\frac{3\alpha_{0}}{2}\right)\right).
\end{split}
\end{equation*}

Equations for the caustics are
\begin{eqnarray*}
x(\rho)&=&x_{0}-\frac{Z_{\pm}^{3}}{3\lambda^{2}A_{0}^{3}\rho^{3m+6}}+\frac{\alpha_{0}Z_{\pm}^{2}}{2\lambda^{2}A_{0}^{2}\rho^{2m+4}}\pm\frac{\alpha_{0}Z_{\pm}}{\lambda\rho}\mp\frac{Z_{\pm}^{2}}{\lambda A_{0}\rho^{m+3}}-\frac{\alpha_{0}A_{0}^{2}\rho^{2m+2}}{2(m+1)},\\
t(\rho)&=&\frac{Z_{+}^{2}}{2\lambda A_{0}^{2}\rho^{2m+4}}+\frac{Z_{+}}{\rho}\left(1-\frac{\alpha_{0}}{\lambda A_{0}\rho^{m+1}}\right)-\alpha_{0}A_{0}\rho^{m+1}+\frac{\lambda A_{0}^{2}\rho^{2m+2}}{2(m+1)}+t_{0},
\end{eqnarray*}
where
\begin{equation*}
Z_{\pm}=\alpha_{0}A_{0}\rho^{m+2}\mp\frac{\lambda A_{0}^{2}\rho^{2m+3}}{2m+3}\pm\alpha_{2}.
\end{equation*}

\begin{figure}[h]

%\sidecaption
\begin{minipage}[h]{0.32\linewidth}
\center{\includegraphics[scale=0.18]{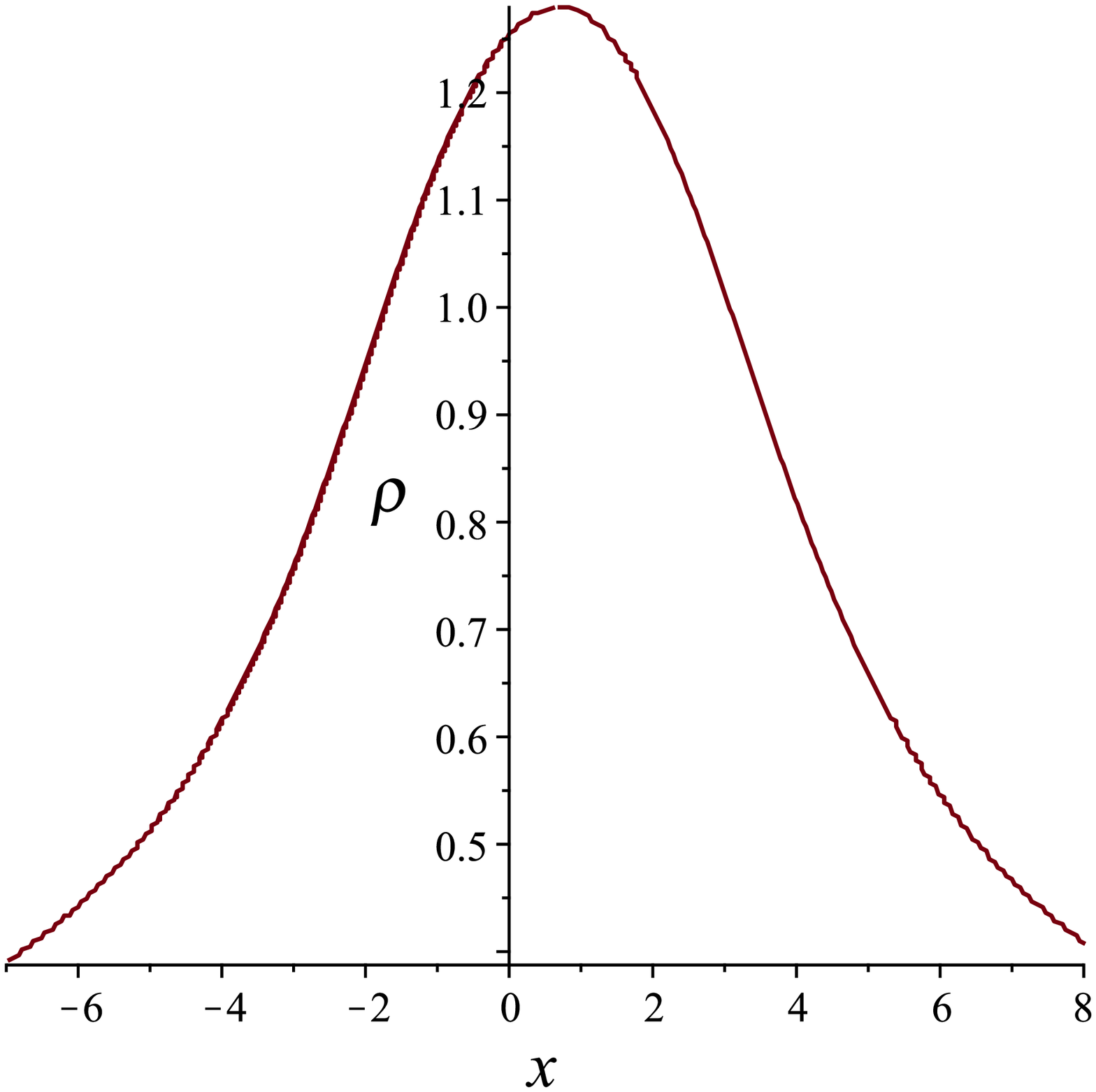}}
\end{minipage}
\hfill
\begin{minipage}[h]{0.32\linewidth}
\center{\includegraphics[scale=0.18]{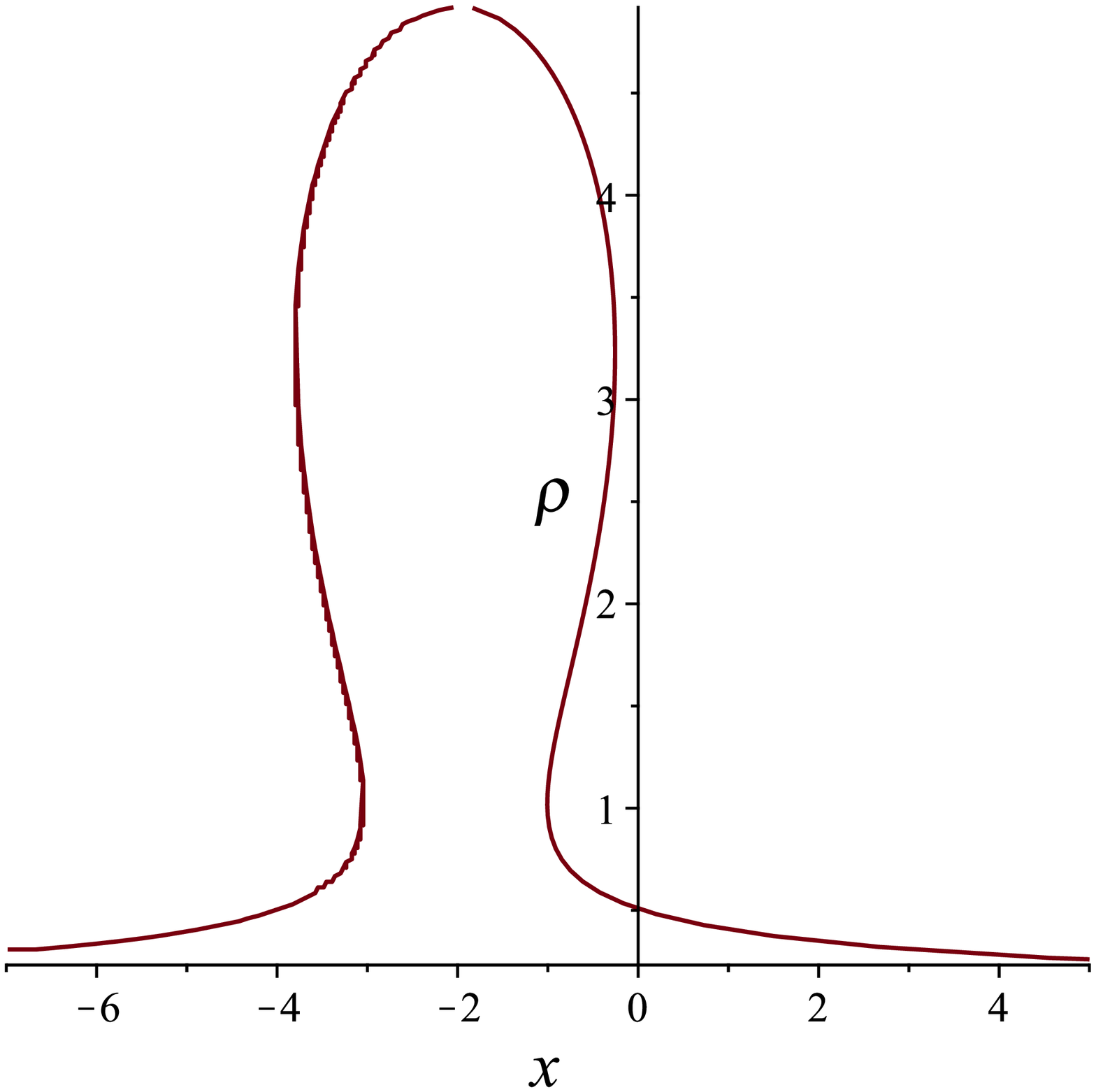}}
\end{minipage}
\hfill
\begin{minipage}[h]{0.32\linewidth}
\center{\includegraphics[scale=0.18]{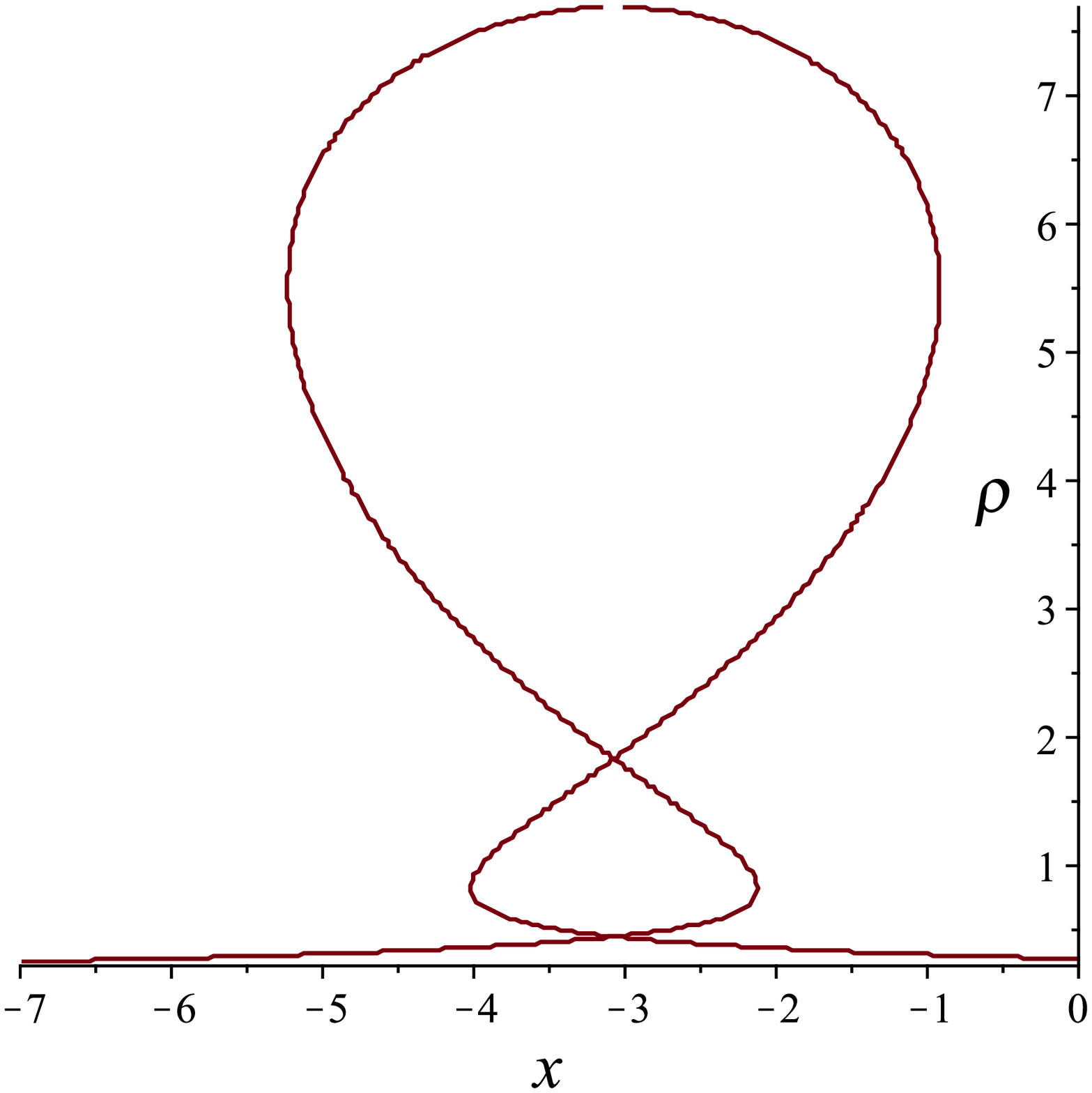}}
\end{minipage}
\caption{Density graph for $n=3$ at moments $t=0$, $t=2.7$, $t=3.75$}
\label{mult}
\end{figure}

The graphs of the multivalued solutions are presented in Fig.~\ref{mult}, where we used substitution $A_{0}=1,\, \alpha_{0}=1,\,t_{0}=1,\,\alpha_{2}=-2,\, x_{0}=0$.

\begin{figure}[h!]
\centering
\includegraphics[scale=.5]{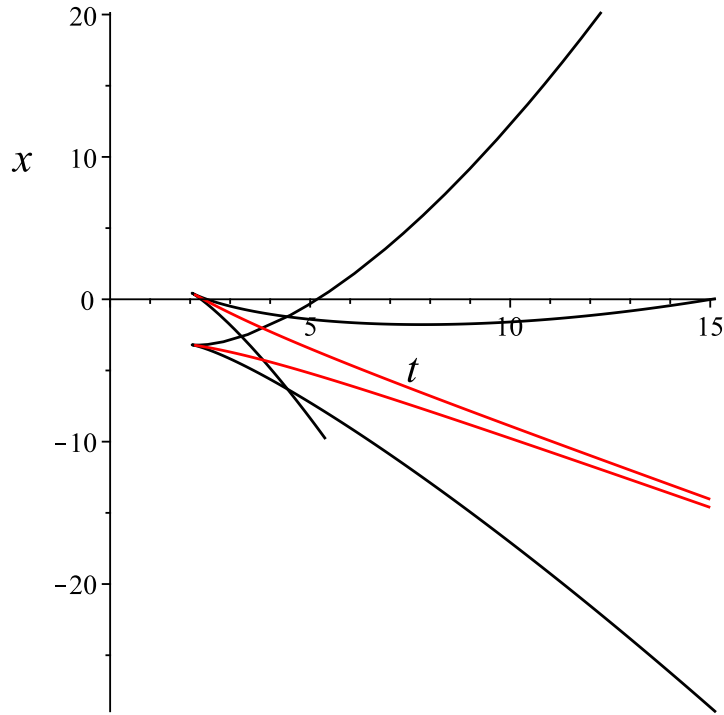}
\caption{Caustics (black) and shock wave fronts (red).}
\label{shock}
\end{figure}

The discontinuity lines together with caustics are shown in Fig.~\ref{shock}. One can observe a very interesting phenomenon. Smooth initial datum evolutes into formation of two cusps, moving towards each other. However, as numerical computations show, two discontinuities meet only asymptotically when $t\to+\infty$.

\section*{Acknowledgements}
The author is grateful to Valentin Lychagin for helpful suggestions during the preparation of the paper. This work was partially supported by the Foundation for the Advancement of Theoretical Physics and Mathematics ``BASIS'' (project 19-7-1-13-3) and by the Russian Science Foundation (project 21-71-20034).

\end{document}